\documentclass[12pt]{article}
\usepackage{amssymb,amsmath,amscd}
\usepackage{xypic}

\newcommand{\cd}{
$$
\spreaddiagramrows{-.7pc}
\spreaddiagramcolumns{-2pc}
\xymatrix{
[a]\in\!\!\!\!\!\! & ~HP_q(X)~ \ar@{^{(}->}[rr] \ar[ddd]^\rho &&  HP_q(X,Z_{\rm
sing}) \\
&& HP_q(\t X) \ulto^\sim_{\pi_*} \dto \rrto && HP_q(\t X,D)
\ulto_{\pi_*}^\sim \dto^{\tilde\rho} & \!\!\!\!\!\!\ni[\t a] \\
&& H^{n-q}(\t X,K_{\t X}) \dlto_{\pi_*}^\sim \rrto &&  H^{n-q}(\t X,K_{\t X}(D)) \\
& H^{n-q}(X,K_X)
}
$$
}

\newcommand{\myfigure}{

\begin{picture}(0,150)(0,0)

\put(110,30){\line(1,0){140}}
\put(230,30){\circle*{3}}
\put(177,20){$\alpha$}
\put(227,15){$D={\rm div}_\infty\alpha$}

\put(230,50){\line(0,1){60}}

\put(225,120){$D_B=D\times{\mathbb C}P^1$}

\qbezier(110,58)(130,60)(150,72)
\qbezier(150,72)(220,110)(250,100)

\put(110,70){\line(1,0){140}}

\put(230,102){\circle*{3}}
\put(230,70){\circle*{3}}

\put(90,75){$A_0$}
\put(90,50){$A_1$}

\put(200,80){\large$\beta$}

\put(137,96){\footnotesize$2\pi i\,{\rm res}_{A_1}\beta=\alpha$}
\put(150,60){\footnotesize$2\pi i\,{\rm res}_{A_0}\beta=-\alpha$}

\put(248,86){\footnotesize
$2\pi i\,{\rm res}_{D_B\cap A_1}({\rm res}_{D_B}\beta)=-{\rm res}_D\alpha$}
\put(248,54){\footnotesize
$2\pi i\,{\rm res}_{D_B\cap A_0}({\rm res}_{D_B}\beta)={\rm res}_D\alpha$}

\put(255,93){\vector(-3,1){20}}
\put(255,61){\vector(-3,1){20}}

\put(0,50){\parbox{20mm}{\Large
$$
\begin{array}{c}
A\times \mathbb CP^1 \\ \downarrow \\ A
\end{array}
$$
}}

\end{picture}}

\newcommand{\mysection}[1]{\section{\large\bf #1}}

\newtheorem{defn}[subsection]{Definition}
\newtheorem{prop}[subsection]{Proposition}
\newtheorem{theo}[subsection]{Theorem}
\newtheorem{lemma}[subsection]{Lemma}

\newenvironment{ssect}[1]{\smallskip\noindent%
\refstepcounter{subsection}%
{\bf \thesubsection~#1}\hspace{-1mm}}
{\smallskip}

\newenvironment{rem}{\smallskip\noindent%
\refstepcounter{subsection}%
{\bf \thesubsection}~~{\sc Remark.}\hspace{-1mm}}
{\smallskip}

\newenvironment{proof}[1]{\noindent {\em Proof#1.}}
{~$\blacksquare$\smallskip}

\newenvironment{ackn}{\medskip \noindent \small
{\sl Acknowledgments.}}{\bigskip}

\newenvironment{smallbibl}[1]
{\small
\begin{thebibliography}{#1}}
{\end{thebibliography}}

\newcommand{\res}{{\rm res}}
\newcommand{\dvsr}{{\rm div}}
\newcommand{\id}{{\rm id}}

\newcommand{\C}{{\mathbb C}}

\newcommand{\OO}{{\cal O}}
\newcommand{\toiso}{\xrightarrow{\sim\;}}

\newcommand{\tpi}{2\pi i}
\newcommand{\dbar}{\bar\partial}
\renewcommand{\t}{\tilde}
\newcommand{\ssc}{sequence of blow-ups with smooth centers}
\newcommand{\ncd}{normal crossing divisor}

\newcommand{\ncsv}{normal crossing subvariety}
\newcommand{\bu}{\mbox{\tiny$_\bullet$}}
\newcommand{\si}{_{\rm sing}}
\newcommand{\WR}{\mbox{\large$\wr$}}
\newcommand{\CP}{\times\C P^1}

\begin{document}


\title{\Large\sc A Polar de Rham Theorem}

\author{Boris Khesin\thanks{
Department of Mathematics,
University of Toronto,
Toronto, ON M5S 3G3, Canada;
e-mail: {\tt khesin@math.toronto.edu}},~
Alexei Rosly\thanks{
Institute of Theoretical and
Experimental Physics,
B.Cheremushkinskaya 25, 117259 Moscow, \mbox{Russia};
e-mail: {\tt rosly@heron.itep.ru}}~
and Richard Thomas\thanks{
Department of Mathematics, Imperial College, London SW7 2AZ, UK;
e-mail: {\tt rpwt@ic.ac.uk}}}

\date{May 5, 2003}

\maketitle

\begin{abstract}
We prove an analogue of the de Rham theorem for polar homology; that the
polar homology $HP_q(X)$ of a smooth projective variety $X$ is
isomorphic to its $H^{n,n-q}$ Dolbeault cohomology group.
This analogue can be regarded as a geometric complexification where arbitrary
(sub)manifolds are replaced by complex (sub)manifolds and de Rham's
operator $d$ is replaced by Dolbeault's $\bar\partial$.
\end{abstract}


\mysection{Introduction} \label{Int}

The idea of polar homology can be explained as follows. In a
complex manifold\footnote{
All manifolds, varieties, their dimension, etc., are understood in this paper
over $\C\,$.}
$X$, consider a $(q+1)$-dimensional submanifold $Y$ and such a
meromorphic $(q+1)$-form $\beta$ on $Y$ that has only first order poles on a
smooth $q$-dimensional submanifold $Z=\dvsr_\infty\beta\subset Y\subset X$.
Under these circumstances, the residue of $\beta$ can be understood as a
holomorphic $q$-form $\alpha=\tpi\,\res\,\beta$ on $Z$ (we include a factor
of $\tpi$ for future convenience).
In other words, to the pair $(Y,\beta)$
we can associate another pair
$(Z,\alpha)=(\dvsr_\infty\beta,\tpi\,\res\,\beta)$ in one dimension less. We
are going to extend this correspondence,
$(Y,\beta)\mapsto(\dvsr_\infty\beta,\tpi\,\res\,\beta)$, to the boundary map
$\partial$ in a certain homological chain complex. Note that if we apply
$\partial$ to the pair $(Z,\alpha)$ above, we get zero because $\alpha$ is
holomorphic. This gives rise to the basic identity $\partial^2=0$. The
formal definition of {\it the polar chain complex} given in the next section
is somewhat lengthier, but its meaning should be already clear. In
particular, the pairs $(Z,\alpha)$ correspond to $q$-cycles if $\alpha$ is a
holomorphic $q$-form on a $q$-dimensional submanifold $Z\subset X$ and such
a cycle is, in fact, a boundary if $\alpha$ is someone's residue.

In the above discussion we considered the situation when only smooth
submanifolds occur. In general, the definition of the polar
chain complex will have
contributions from arbitrary subvarieties $Z\subset X$. Such a definition,
which gives us a chain complex with homology groups to be denoted as $HP_q(X)$,
was suggested in refs. \cite{KR,KR2}. In many aspects it is analogous to the
definition of topological homology (say, singular homology). In the present
paper, we are going to prove a theorem analogous to de Rham's theorem in the
topological context. Namely, we shall prove that the groups $HP_q(X)$
for smooth projective $X$ are dual to $H^q(X,\OO_X)$, as it was conjectured
in ref. \cite{KR}. In other words, we shall see that the Dolbeault
$\dbar$-complex on $(0,q)$-forms interacts with the polar chain complex in
the same way as the de Rham $d$-complex does with ordinary topological
chains.
The reader interested only in reading the main results should, after having
a look at the definition \ref{PHdef}, proceed directly to Theorem \ref{deRh}
and its proof in \ref{proof}. The rest of the paper consists of technical
preliminaries needed to deal with singularities.

One should note that there exists a more general polar complex, where the
chains are complex subvarieties of dimension $q$ with logarithmic
$p$-forms on them. The corresponding polar homology groups, enumerated by
two indices, are, in general, not isomorphic to any Dolbeault
homology as simple examples show. From this point of view, the isomorphism for $p=q$
discussed in this paper is rather an exception than a rule.

The motivation for considering polar homology comes
from mathematical physics. It appears naturally in ``holomorphization" of
various topological objects; cf. \cite{DT, KR2}.


\mysection{Definitions}\label{Def}

\begin{ssect}{Poincar\'e residue.}
Let $X$ be a smooth complex projective $n$-dimensional manifold and
$V\subset X$ a smooth hypersurface in $X$. Consider a meromorphic $n$-form
$\omega$ on $X$ with first order poles on $V$. If $\{z=0\}$ is a local
equation for $V$, the form $\omega$ can be written as $$ \omega =
\frac{dz}{z}\wedge\rho + \gamma \,, $$ where the locally defined holomorphic
forms $\rho$ and $\gamma$ can be chosen in various ways. However, the
restriction of $\rho$ to $V$ is defined uniquely and, therefore, becomes a
global holomorphic $(n-1)$-form on $V$. It is denoted by
$\res\,\omega=\rho|_V$ and is called the {\it Poincar\'e residue} of
$\omega$. This can be also described by the following exact sequence of
sheaves:
\begin{equation} \label{sh.seq.}
0 \to K_X \to K_X(V) \to K_V \to 0  \,,
\end{equation}
where $K_X$ is the canonical sheaf on $X$, i.e., the sheaf of holomorphic
$n$-forms, while $K_X(V)$ stands for  $n$-forms with first order poles on
$V$ whose residues give us regular $(n-1)$-forms on $V$. The restriction map
$K_X(V) \to K_V$ represents here the Poincar\'e residue for locally defined
$n$-forms. The corresponding residue map for the globally defined forms,
$res: H^0(X,K_X(V)) \to H^0(V,K_V)$, shows up in the cohomological long
exact sequence implied by (\ref{sh.seq.}):
\begin{equation} \label{l.seq.}
\begin{array}{ll}
0 \to & H^0(X,K_X) \to H^0(X,K_X(V))
      \stackrel{\!\!res}{\longrightarrow} H^0(V,K_V) \to  \\
   &   H^1(X,K_X) \to H^1(X,K_X(V)) \longrightarrow \dots
\end{array}
\end{equation}
In this sequence we encounter elements of polar homology. Namely, the meromorphic
$n$-forms $\omega\in H^0(X,K_X(V))$ will correspond
(via the definitions in \ref{PHdef} below)
to $n$-chains, the holomorphic $(n-1)$-forms $\rho\in H^0(V,K_V)$ will
correspond to $(n-1)$-cycles, while the boundary map will be given by the
map $res$ in (\ref{l.seq.}). We shall see that the contribution to the
$(n-1)$-dimensional polar homology coming from a given
(smooth) hypersurface $V$ will correspond to the quotient
$H^0(V,K_V)/res(H^0(X,K_X(V)))$. It remains to understand the
contributions from arbitrary subvarieties in $X$.

\end{ssect}

\begin {ssect}{Normal crossings.}
Since we are going to use the map $res: H^0(X,K_X(V)) \to H^0(V,K_V)$ in the
definition of a boundary map on a vector space of chains we cannot restrict
to the case of only smooth divisors of poles. As a matter of fact, it is
sufficient to generalize to the case of normal crossings. We shall consider
\ncd s, as well as subvarieties with normal crossings in arbitrary
codimension.
We shall give a very restrictive definition of these which will suffice
for our purposes.
Let us explain our conventions in more detail. First of all, a (sub)variety
will be always reduced, but not necessary irreducible.
Thus, a subvariety\footnote{ In this paper the varieties are always
projective or quasi-projective; the subvarieties are always closed.} in $X$
is just a Zariski closed subset of $X$. On the other hand, a smooth variety
(= smooth manifold = manifold) will be always assumed irreducible
(which is equivalent to connected for smooth varieties).

Let us consider a smooth $n$-dimensional manifold $X$. A hypersurface
$V\subset X$ will be called a {\it normal crossing divisor} if $V$ consists
of smooth components that meet transversely, in the sense that $V=\cup_iV_i$,
where each $V_i$ is smooth and intersects transversely $V_j$, $V_j\cap V_k$, and
so on, for all $i,j,k,\ldots$.\footnote{
Near each point $x\in V$, one can choose local coordinates
$z_1,\dots,z_n$ in $X$ in such a way that $z_1\cdot \ldots \cdot z_p=0$ is a
local equation of $V$ (where $p\leqslant n$ is the number of components of $V$ passing
through $x$). The latter local formulation could be used as a definition of a
\ncd. We prefer, however, a stronger version, when the self-intersections of
components are excluded.}
In order to introduce the notion of a \ncsv\ of an arbitrary codimension,
consider first a codimension two subvariety $W\subset V\subset X$ (where
$X$ and $V$ are as above). Let us require that the part of $W$ which resides
in a smooth component of $V$ is a \ncd\ there and that $W$ intersects the normal
crossing singularities of $V$ transversely.
More precisely, if $W\nsupseteq V_i\cap V_j, \forall i,j$, and
$(W\cap V_i)\cup(V_i\cap(\cup_{k\neq i}V_k))$
is a \ncd\ in the smooth manifold $V_i$ for all $i$, we shall say that $W$
is a \ncd\ in $V$ and a \ncsv\ in $X$.
In such a way we obtain
the notion of a normal crossing divisor in a variety, which
is itself a normal crossing divisor in a bigger variety.
Proceeding deeper in codimension we shall say that a subvariety $Y$ of
codimension $m$ in $X$ is a \ncsv\ if there exists a nested sequence
\begin{equation} \label{flag}
Y=V^m\subset V^{m-1}\subset\ldots\subset V^1\subset V^0=X \,,
\end{equation}
such that $V^{i+1}$ is a \ncd\ in $V^i\,$.
We shall also say that two \ncd s $V$ and $V'$ intersect transversely if
$V+V'$ is a \ncd\ again.
(This means in particular that $V$ and $V'$ have no common components and
that $V\cap V'$ is a \ncd\ both in $V$ and in $V'$.)

In fact, we shall need mainly the notion of an {\it ample} subvariety
with normal crossings in a projective manifold $X$.

\end{ssect}
\begin{defn} \label{ANCV}
A \ncsv\ $Y\subset X$ in a projective manifold $X$ is called ample if one
can choose a flag (\ref{flag}) in such a way that $V^{i+1}$ is an ample
\ncd\ in $V^i\,$.

\end{defn}
\begin {ssect}{Canonical line bundle.}
The canonical sheaf $K_V$ is defined for a smooth variety $V$ as the sheaf
of holomorphic forms of the top degree on $V$ and, if $V$ is a hypersurface
in some $X$, $i\!:V\hookrightarrow X$, the local properties are described by
the sequence (\ref{sh.seq.}). In this case, one {\it has to show} that
$i^*K_X(V)\simeq K_V$, while the Poincar\'e residue gives us a canonical
choice of this isomorphism. In the case of a normal crossing divisor
$i\!:V\hookrightarrow X$ we may take the sequence (\ref{sh.seq.}) as the
{\it definition} of $K_V$. In other words, $K_V$ is defined as $i^*K_X(V)$.
By induction in codimension we obtain a definition that can be applied to
any \ncsv\ $Y$; the result is a line bundle on $Y$ which does not depend on
the choice of the flag (\ref{flag}): invariantly,
$K_Y={\cal E}xt^m(\OO_Y,K_X)$, where $m={\rm codim}\,Y$.
With such a definition, the global
sections of $K_V$ are regarded as ``holomorphic" forms on $V$ and the
Poincar\'e residue, $res: H^0(X,K_X(V)) \to H^0(V,K_V)$, still maps
meromorphic forms to holomorphic ones. This is precisely what we need to
define a chain complex.

As a last preparation, it remains to check only the properties of the
repeated residue map, as it has to support the identity $\partial^2=0$. Let
$V$ be a \ncd\ and suppose for simplicity that it consists of only two
components, $V=V_1\cup V_2$, so that $V_1,V_2$ are smooth and intersect
transversely over a smooth variety $V_{12}=V_1\cap V_2$. Then, a section
$\alpha\in K_V$ can be described via its restrictions
$\alpha_i=\alpha|_{V_i}$. Since $K_V|_{V_i}\simeq K_X(V_1+V_2)|_{V_i}\simeq
K_X(V_i)|_{V_i}(V_1\cap V_2)\simeq K_{V_i}(V_{12})$, the $\alpha_i$ are in
fact meromorphic forms, $\alpha_i\in H^0(K_{V_i}(V_{12}))$. Moreover, it follows
from a local coordinate calculation with the definition that $\res_{V_{12}}
\alpha_1+\res_{V_{12}}\alpha_2=0$,
which is summarized in the short exact sequence of sheaves
$$
0 \to K_V \to K_{V_1}(V_{12}) \oplus K_{V_2}(V_{12}) \to K_{V_{12}} \to 0\,,
$$
where the third arrow is taking the sum of residues. In other words, a
holomorphic form $\alpha\in H^0(V,K_V)$ on a normal crossing variety $V$ can
be described as a collection of meromorphic forms $\alpha_i$ on $V_i$
satisfying the pairwise cancellation of their residues at the intersections.
(We shall say that the polar cycle $(V,\alpha)$ is the sum of two polar
chains $(V_1,\alpha_1)$ and  $(V_2,\alpha_2)$, whose boundaries cancel each
other.)
\end{ssect}

\begin{ssect}{Resolution of singularities.} \label{Hironaka}
In the next section, our main tool will be the Hironaka
theorem on resolution of singularities \cite{Hi}. This theorem asserts that
every algebraic variety $Z$ admits a desingularization, that is there exists
a smooth variety $\t Z$ and a regular projective birational morphism
$\pi:\t Z\to Z$, which is biregular over $Z-Z\si\,$. Moreover, $\pi$ can be
obtained as a sequence of blowing up with smooth centers. If $D$ is a
subvariety in $Z$ we can additionally require that $\pi^{-1}(D)$ is a \ncd\
in $\t Z$.

We shall also need the following important result, the (weak)
factorization theorem for birational morphisms, proved recently by
Abramovich, Karu, Matsuki and W\l odarczyk \cite{W,A...W}. Below we cite only a part
of their statement from ref. \cite{A...W} relevant to our needs
(the complete proposition is much stronger).
\end{ssect}

\begin{prop} \label{A...W}
Let $\phi:X\dashrightarrow X'$ be a birational map between smooth projective
varieties $X$ and $X'$. Then $\phi$ can be factored into a sequence of
blowings up and blowings down with smooth irreducible centers, namely, there
exists a sequence of birational maps between smooth projective varieties
$$
X=\t X_0  \stackrel{\varphi_1}{\dashrightarrow}
\t X_1 \stackrel{\varphi_2}{\dashrightarrow}
\cdots \stackrel{\varphi_i}{\dashrightarrow}
\t X_i \stackrel{\varphi}{\dashrightarrow}
\t X_{i+1} \stackrel{\varphi_{i+2}}{\dashrightarrow}
\cdots \stackrel{\varphi_{l-1}}{\dashrightarrow}
\t X_{l-1} \stackrel{\varphi_l}{\dashrightarrow}
\t X_l = X'
$$
where
\begin{enumerate}
\item
$\phi=\varphi_l\,\circ\varphi_{l-1}\,\circ\cdots\circ\varphi_2\,\circ\varphi_1\,$,
and

\item
either
$\varphi_i:\t X_{i-1}\dashrightarrow\t X_i$, or
$\varphi_i^{-1}:\t X_i\dashrightarrow\t X_{i-1}$
is a morphism obtained by blowing up a smooth irreducible center.

\end{enumerate}
\end{prop}

For the sake of brevity in what follows, under a `blow-up' we shall
understand `a sequence of blowings up with smooth centers'.
The following corollary of the Hironaka and Bertini theorems will also be useful
in the sequel.

\begin{prop} \label{Hi-flag}
Let $Z\subset X$ be an arbitrary irreducible subvariety of codimension $m$
in a smooth projective manifold $X$. Then, there exists a blow-up
$\pi:\t X\to X$ and a flag of subvarieties
\begin{equation} \label{Z-flag}
\t Z\subset V^{m-1}\subset V^{m-2}\subset\ldots\subset V^1\subset V^0=\t X
\end{equation}
such that $V^{i+1}$ is a smooth hypersurface in $V^i$ and
$\t Z$ is smooth and mapped birationally by $\pi$ onto $Z$.
\end{prop}

\begin{proof}{}
Firstly, by Hironaka, we can blow up $X$ in such a way that the proper
preimage of $Z$ becomes smooth. If the codimension of $Z$ is one, $m=1$, the
proposition is proved. We can thus proceed for $m>1$ and assume that $Z$ is
already smooth. In this case, let us take a very ample divisor class $H$ in
$X$ and consider hypersurfaces in this class containing $Z$. Such
hypersurfaces are described as zero sets of global sections of the sheaf
${\cal I}_Z(H)$, where ${\cal I}_Z\subset \OO_X$ is the ideal sheaf of the
subvariety $Z$ in $X$. By Bertini, the generic section
$s\in H^0(X,{\cal I}_Z(H))$ defines a hypersurface $V=\{s=0\}\subset X$ which is
regular outside $Z$. As to the points of $V$ which lie on $Z$, the
singularities correspond to the zeros of the section
$\bar s=ds\in H^0(Z,{\cal I}_Z/{\cal I}_Z^2(H))$ induced by $s$. Let us choose
$H\gg0$ in such a way that $H^0(Z,{\cal I}_Z/{\cal I}_Z^2(H))\neq 0$, while
$H^1(Z,{\cal I}_Z^2(H))=0$. Then we have a non-trivial section
$\bar s$ in $H^0(Z,{\cal I}_Z/{\cal I}_Z^2(H))$ whose zeros form a proper
closed subset $Z_0\varsubsetneq Z$. Moreover,
$H^1(Z,{\cal I}_Z^2(H))=0$ guarantees that the mapping
$H^0(X,{\cal I}_Z(H))\to H^0(Z,{\cal I}_Z/{\cal I}_Z^2(H))$,
$s\mapsto\bar s$, is surjective. Hence, taking a generic $s$ we can ensure
that the resulting hypersurface $V=\{s=0\}$ is regular outside
$Z_0=\{\bar s=0\}\varsubsetneq Z$. Applying the Hironaka theorem, we can now
resolve the singularities of $V$ by blowing up $X$ in centers belonging to
$Z_0\subset X$. Then, for the proper preimage $\t Z$ of $Z$, we have that
$\t Z\subset V^1\subset\t X$, where $\t Z$ and $V^1$ are smooth.
We can then proceed in the same manner inside $V^1$ until
the whole flag (\ref{Z-flag}) obeying the required conditions is
constructed.
\end{proof}

\begin{ssect}{Polar chains.}
The space of polar $q$-chains for a (not necessarily smooth) complex projective
variety $X$, $\dim X=n,$ will be defined as a $\C\,$-vector space with
certain generators and relations.

\end{ssect}

\begin{defn} \label{PHdef}
The space of polar $q$-chains~ ${\cal C}_q(X)$ is
a vector space over $\C$ defined as the quotient
${\cal C}_q(X)=\hat{\cal C}_q(X)/{\cal R}_q$, where the vector space
$\hat{\cal C}_q(X)$ is freely generated by the
triples $(A,f,\alpha)$ described in (i),(ii),(iii) and
${\cal R}_q$ is defined as relations
{\rm (R1),(R2),(R3)} imposed on the triples.
\begin{description}
\item (i) $A$ is a smooth complex projective variety, $\dim A=q$;
\item (ii) $f\!: A\to X$ is a holomorphic map of projective varieties;
\item (iii) $\alpha$ is a meromorphic $q$-form on $A$
with first order poles on $V\subset A$, i.e., $\alpha\in H^0(A,K_A(V))$,
where $V$ is a normal crossing divisor in $A$.
\end{description}

\noindent
The relations are generated by:
\begin{description}
\item {\rm (R1)} $\lambda (A, f, \alpha)=(A, f, \lambda\alpha)$,

\item {\rm (R2)}
$\sum_k(A_k,f_k,\alpha_k)=0$ provided that
$\sum_kf_{k*}\alpha_k\equiv 0$ on a Zariski open dense subset of $\hat A$,\footnote{
For a surjective holomorphic map $f:U\to V$ of two smooth complex manifolds
of the same dimensions (that is to say, $f$ is generically finite), we have
a push-forward map $f_*$ on differential forms defined on the locus over
which $f$ is finite by the summation over
the preimages $P\in f^{-1}(Q)$ of a point $Q$ . This map is also called
the trace map, and the pushforward of holomorphic (resp. meromorphic) forms
extend over the image to be holomorphic (resp. meromorphic) \cite{Griff}.}
where $f_k(A_k)=f_l(A_l)=:\hat A, ~\forall\, k,l$ and
$\dim\hat A=\dim f_k(A_k)=q, ~\forall k$;

\item {\rm (R3)} $(A,f,\alpha)=0$ if $\dim f(A)<q$.
\end{description}

\end{defn}

\begin{defn} \label{d-def}
The boundary operator~ $\partial: {\cal C}_q (X)\to{\cal C}_{q-1}(X)$
is defined by
$$
\partial(A,f,\alpha)=\tpi\sum_k(V_k, f_k, \res_{V_k}\,\alpha) \,,
$$
where $V_k$ are the
components of the polar divisor of $\alpha$,
$\dvsr_\infty\alpha=\cup_kV_k$, and the maps $f_k=f|_{V_k}$
are restrictions of the map $f$ to each component of the divisor.

\end{defn}

\begin{prop}
The boundary operator $\partial$ is well defined, i.e.\ it is
compatible with the relations {\rm (R1),(R2),(R3)}.

\end{prop}

\noindent
For the proof see \cite{KR}. Now, by using the cancellation of repeated
residues for  forms $\alpha$ with normal crossing divisors of poles, one
proves the following \cite{KR}:

\begin{prop}
~~~$\partial^2=0\;$.
\end{prop}

\noindent
This allows one to define a homology theory.

\begin{defn}
For a complex projective variety $X, \dim X=n$, the chain
complex
$$
0 \to {\cal C}_n(X) \xrightarrow{\;\partial~}
{\cal C}_{n-1}(X) \xrightarrow{\;\partial~} \dots
\xrightarrow{\;\partial~} {\cal C}_0(X)\to 0
$$
is called the polar chain complex of $X$. Its homology groups,
$HP_q(X), q=0,\ldots,n$, are called the polar homology groups of $X$.

\end{defn}

\begin{rem}\label{supp}
It is useful to introduce the notion of the support of a $q$-chain
$a\in{\cal C}_q(X)$. This is defined as the following minimal subvariety
${\rm supp}\,a=\bigcap\cup_kf_k(A_k)\subset X$ where the intersection
$\bigcap$ is taken over all representatives $\sum_k(A_k,f_k,\alpha_k)$ in
the equivalence class $a$. (In other words, ${\rm supp}\,a$ can be
determined by taking $Z=\cup_kf_k(A_k)$ for an arbitrary representative
$\sum_k(A_k,f_k,\alpha_k)$, removing those components of $Z$ which are of
dimension less than $q$ or where the push-forwards $f_{k*}\alpha_k$ sum to
zero as in (R2) in the Definition \ref{PHdef} above and taking closure.)
This notion of the support of a polar chain coincides with the support of
the current in $X$ corresponding to that chain. (The relation with currents
was discussed in ref.\ \cite{KR}.)

If $a\in{\cal C}_q(X)$ then $Z={\rm supp}\,a$ is either of pure dimension
$q$, or empty. The smooth part of $Z$ is provided with a meromorphic $q$-form
$\alpha$ obtained by summation of $f_{k*}\alpha_k$. The meaning of the relation
(R2) above is essentially that these data, $({\rm supp}\,a,\alpha)$, define the
equivalence class of (sums of) triples $a\in{\cal C}_q(X)$ in a unique way.
By the Hironaka theorem, the subvariety $Z$ can in fact be arbitrary, that is for
an arbitrary $q$-dimensional $Z\subset X$, there exists a $q$-chain $a$ such that
$Z={\rm supp}\,a$, but the meromorphic $q$-form $\alpha$ on $Z-Z_{\rm sing}$ cannot
in general be arbitrary.

\end{rem}
\begin{ssect}{Relative polar homology.} \label{relPHdef}
Let $Z$ be a closed subvariety in a projective $X$.
Analogously to the topological relative homology we can define the
polar relative homology of the pair $Z\subset X$.

\end{ssect}

\begin{defn} The relative polar homology groups $HP_q(X,Z)$ are the
homology groups of the following quotient complex of chains:
$$
{\cal C}_q(X,Z)\:={\cal C}_q(X)/{\cal C}_q(Z).
$$
\end{defn}

\noindent
Here we use the natural embedding of the chain groups
${\cal C}_q(Z)\hookrightarrow{\cal C}_q(X)$.
This leads to the long exact sequence in polar homology:
\begin{equation} \label{l.seq.rel.}
\ldots\to HP_{q}(Z)\to HP_q(X)\to HP_q(X,Z)
\xrightarrow{\;\partial~}
HP_{q-1}(Z)\to\ldots
\end{equation}

\begin{ssect}{}
The functorial properties of polar homology are straightforward.
A regular morphism of projective varieties
$h\!:X\to Y$ defines a homomorphism
$h_*\!:HP\bu(X)\to HP\bu(Y)$.
Analogously, for the relative polar homology we have
$h_*\!:HP\bu(X,V)\to HP\bu(Y,W)$ if $V\subset X$, $W\subset Y$ are closed
subsets and $h(V)\subset W$.

\end{ssect}

\begin{rem} \label{2pairs}
In the case of a morphism of two pairs $h\!:(X,V)\to(X',V')$ as above, the
induced homomorphisms $h_*$ give us the homomorphism of the associated long
exact sequences:
\begin{equation}
\begin{array}{ccccccccc}
\ldots\to& HP_q(V)\, &\to& HP_q(X)\, &\to& HP_q(X,V)\,
&\to& HP_{q-1}(V)\, &\to\ldots \\
& \downarrow && \downarrow && \downarrow && \downarrow & \\
\ldots\to& HP_q(V')  &\to& HP_q(X') &\to& HP_q(X',V')
&\to& HP_{q-1}(V') &\to\ldots
\end{array}
\end{equation}
We note that if any two of the three homomorphisms
$HP\bu(V)\to HP\bu(V')$, $HP\bu(X)\to HP\bu(X')$,
$HP\bu(X,V)\to HP\bu(X',V')$ are isomorphisms then
the third one is an isomorphism as well.

\end{rem}


\mysection{Polar Homology and Dolbeault cohomology}
\label{Hom}

We are going to show that the Dolbeault, or $\dbar$, cohomology on
$(0,q)$-forms, $H^{(0,q)}_{\dbar}(X)$, plays the same role with respect to
polar homology $HP_q(X)$ as does the de Rham cohomology in the topological
context. First of all, there is an obvious pairing between $HP_q(X)$ and
$H^{(0,q)}_{\dbar}(X)$. For $[(A,f,\alpha)]\in HP_q(X)$ and $[\omega]\in
H^{(0,q)}_{\dbar}(X)$, we can write $\int_A\alpha\wedge f^*\omega$ and show
that such a pairing descends to (co)homology classes. Recalling the
isomorphism $H^{(0,q)}_{\dbar}(X)\simeq H^q(X,\OO_X)$ and by the Serre
duality, $H^q(X,\OO_X)^*\simeq H^{n-q}(X,K_X)$, the above pairing is thus
represented by the map
\begin{equation}\label{rho}
\rho : HP_{q}(X)\to H^{n-q}(X, K_X) \,,
\end{equation}
where $n=\dim X$.

\begin{theo}\label{deRh} {\bf (Polar de Rham theorem)}
For a smooth projective $n$-dimensional $X$,
the map $\rho$ is an isomorphism for any $q$:
$$
HP_{q}(X)\simeq H^{n-q}(X, K_X) \;.
$$
\end{theo}

\noindent
In the case of polar homology of $X$ relative to a hypersurface $V\subset X$
we analogously have the pairing of $HP_{q}(X,V)$ and $H^q(X,\OO_X(-V))$, or,
by Serre's duality, the homomorphism
\begin{equation}\label{rho rel}
\rho : HP_{q}(X,V)\to H^{n-q}(X, K_X(V)) \,,
\end{equation}
and the corresponding relative version of the Theorem \ref{deRh} is as
follows.

\begin{theo}\label{rel} Let $V$ be a normal crossing divisor in a smooth
projective $X$. Then
$$
HP_{q}(X,V)\simeq H^{n-q}(X, K_X(V)) \;.
$$
\end{theo}

\noindent This more general assertion follows in fact from the Theorem
\ref{deRh} by comparing the long exact sequence in sheaf cohomology
(\ref{l.seq.}) with that in relative polar homology, cf.,
(\ref{l.seq.rel.}).

\begin{rem}
It follows from Theorem \ref{deRh} that if two smooth projective manifolds $X$
and $X'$ are birationally equivalent, then $HP_q(X)=HP_q(X')$ since we have in
this case that $H^{n-q}(X, K_X)=H^{n-q}(X', K_{X'})$. However, we in fact
prove this and other similar results first without
reference to sheaf cohomology, on the way to the proof of Theorem \ref{deRh}.
In fact the rest of the paper is now devoted to proving Theorem \ref{deRh}.
\end{rem}

\begin{lemma}\label{BRE}
If two projective varieties $X$ and $X'$ are birationally equivalent and we
have an isomorphism
$$
g: X-Z \toiso X'-Z' \,,
$$
where $Z$ (resp.\ $Z'$) is a Zariski closed subset in $X$ (resp.\ in $X'$),
then
$$
HP\bu(X,Z) \simeq HP\bu(X',Z') \,.
$$
\end{lemma}

\begin{proof}{}
We want to construct an isomorphism of complexes
\begin{equation} \label{g}
g\bu:{\cal C}\bu(X,Z) \toiso {\cal C}\bu(X',Z').
\end{equation}
Let us take an arbitrary non-zero simple\footnote{
We call a chain simple if it is equivalent to a single triple rather than a
sum of triples.}
chain $a\in {\cal C}_q(X,Z)$ and let the triple $(A,f,\alpha)$ be a
representative of the equivalence class $a$. Since $a\neq0$, the
image $\hat A=f(A)$ of $A$ in $X$ has $\dim\hat A=q$ and
$\hat A \nsubseteq Z$. Let us define $\hat A'$ as the closure of
$g(\hat A-Z)$ in $X'$. By the Hironaka theorem (take the closure of
the graph of $g|_{A-Z}$ in $A\times \hat A'$ and resolve),
there exists a smooth $q$-dimensional variety $A'$ with regular maps
$f':A'\to X'$ and $\pi:A'\to A$, where $\pi$ is a birational map of $A'$
onto $A$, such that they form together with $f$ and $g$  (on open dense
subsets) a commutative square, namely:
\begin{equation} \label{sq}
\begin{CD}
A-f^{-1}(Z)    @>f>>   \hat A-Z  @. ~\hookrightarrow~ @. X-Z   \\
  @A\pi AA              @VV\WR V    @.       @VgV\WR V               \\
A'-{f'}^{-1}(Z)  @>{f'}>>   \hat A'-Z' @. ~\hookrightarrow~ @. X'-Z'
\end{CD}
\end{equation}
By blowing up $A'$ further if necessary and setting $\alpha':=\pi^*\alpha$, we
may assume that $\dvsr_\infty\alpha'$ is a \ncd, along which $\alpha'$ has
\emph{first order} poles. This is because it is a top degree form, for which
having first
order poles is the same as being logarithmic, and logarithmic forms are
locally generated as a ring by forms $df/f=d\log f$ which are also logarithmic
on pullback. So $(A',f',\alpha')$ is admissible and
defines a chain $a'\in {\cal C}_q(X',Z')$. We define the map $g_q$
by setting $a'=g_q(a)$.

Note that the $q$-forms $f_*\alpha$ and $f'_*\alpha'$, which are defined on
open dense subsets in $\hat A$ and $\hat A'$ respectively, coincide there
(in the sense of the isomorphism $g:\hat A-Z\toiso\hat A'-Z'$) as follows
from the commutative diagram (\ref{sq}). This observation shows us that
$g_q:a\mapsto a'$ is well defined, because, in general, polar chains are
uniquely defined in terms of the forms $f_*\alpha$ on the dense subsets in
their supports (cf.\ Remark \ref{supp}). It is obvious that the same
construction applied to $g^{-1}:X'-Z'\toiso X-Z$ gives the inverse of
$g\bu\,$. Compatibility with the boundary map $\partial$ is also obvious. Thus
we have indeed constructed an isomorphism of complexes (\ref{g}), which
proves the lemma.
\end{proof}
\begin{lemma}\label{cyl}
Let $M$ be any projective variety, then
$$
HP\bu(M\CP) \simeq HP\bu(M) \,,
$$
where the isomorphism is induced by the projection $\pi:M\CP\to M$.
\end{lemma}

\begin{proof}{}\
Choosing a point $0\in\C P^1$, we will show that any cycle
in $M\CP$ is homologous to one in the zero section $s=(\id,0):
M\to M\CP$ by constructing a homotopy
$h:{\cal C}_q(M\CP)\to{\cal C}_{q+1}(M\CP)$ from $s_*\circ\pi_*$ to the
identity; that is
\begin{equation} \label{dh+hd}
\partial\circ h + h\circ\partial = \id - s_*\circ\pi_* \,.
\end{equation}

Let $a=(A,f,\alpha)\in{\cal C}_q(M\CP)$ be a simple chain; that is $\dim A=q$,
$\alpha$ is a $q$-form on $A$ whose poles form a \ncd\ in $A$, and $f=(f_M,g)$
with $f_M:=\pi\circ f:A\to M$ a regular map and $g:A\to\C P^1$ a rational
function on $A$. We would like to define the $(q+1)$-chain $h(a)$ by
$$
h(a)=(A\CP,f_M\times{\rm id}_{\C P^1},\beta), \qquad\text{where}\quad
\beta= \frac{1}{\tpi}\, \frac{g\,dz}{z(z-g)}\wedge\alpha \,.
$$
Here $z$ is an inhomogeneous coordinate on $\C P^1$ vanishing at
$0\in\C P^1$, and $z,g$ and $\alpha$ are pulled back to the product
$A\CP$. $\beta$ has simple
poles on the hypersurface
$\dvsr_\infty\beta=A_1\cup A_0\cup(\dvsr_\infty\alpha\CP)$, where
$A_1=\{z=g\}$ and $A_0=\{z=0\}$ are two sections, so that, in particular,
$A_1\simeq A_0\simeq A$.

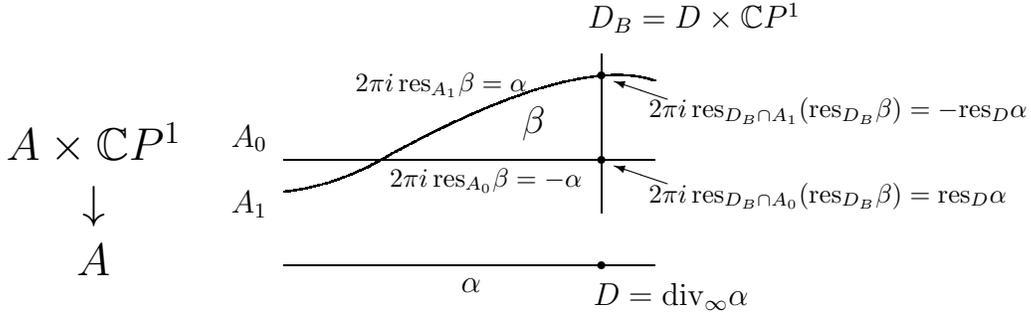
\myfigure

\noindent
The corresponding residues are as follows:
\begin{equation} \label{resb}
\begin{array}{ccl}
\tpi\,\res_{A_1}\,\beta & = & ~\alpha \,, \\
\tpi\,\res_{A_0}\,\beta & = & -\alpha \,, \\
\tpi\,\res_{\dvsr_\infty\alpha\CP}\,\beta & = &
-\frac{g\,dz}{z(z-g)}\wedge\res\,\alpha \,.
\end{array}
\end{equation}
The only problem is that $\dvsr_\infty\,\beta$ will not be a \ncd\ if $A_0$
does not meet $A_1$ or $A_1\cap(\dvsr_\infty\alpha\CP)$ transversely.

By changing $z$ to $z'$ (and
so moving $0\in\C P^1$) we can ensure that the new $A_0'$ \emph{does} meet
$A_1$ and $A_1\cap(\dvsr_\infty\alpha\CP)$ transversely, and the resulting
$\beta'$ has normal crossing poles, but now the definition of $h'(a)$ appears
to depend on the choice of $A'_0$.
The solution is to take
this new $\beta'$ and add to it $\beta-\beta'$, which also has normal crossing
poles (along $A_0\cup A_0'\cup(\dvsr_\infty\alpha\CP)$). Thus
$h(a)=\beta=\beta'+(\beta-\beta')$
is admissible in the sense of Definition \ref{PHdef}, and $h$ is well defined
and linear.

>From (\ref{resb}) we can now calculate:
\begin{equation} \label{dh}
\begin{array}{ccl}
\partial h(a) & = &
 (A_1\,, (f_M\times{\rm id}_{\C P^1})|_{A_1}\,,\alpha) \,-\,
 (A_0\,, (f_M\times{\rm id}_{\C P^1})|_{A_0}\,,\alpha) \\
&& \hspace{15mm}-\ (\dvsr_\infty\,\alpha\CP,
f_M|_{\dvsr_\infty\,\alpha}\times{\rm id}_{\C P^1}\,,
\frac{g\,dz}{z(z-g)}\wedge\res\,\alpha) \\
&=& (A\,,f\,,\alpha)\,-\,(A\,,s\circ\pi\circ f\,,\alpha)\,-\,
h(\dvsr_\infty\,\alpha,
f|_{\dvsr_\infty\,\alpha}\,,\tpi\,\res\,\alpha) \\
&=& a\,-\,s_*\pi_*(a)\,-\,h\partial(a)\,,
\end{array}
\end{equation}
as in (\ref{dh+hd}).
\end{proof}

\begin{lemma}\label{bl-up}
~\\
(a) Let $M$ be a smooth projective variety and $E$ be the total space of
a projective bundle over $M$, i.e.\ $\pi:E\to M$ is a locally trivial
fibration (in the Zariski topology) with a projective space as a fiber. Then
$\pi$ induces an isomorphism in polar homology:
$$
HP\bu(E) \simeq HP\bu(M) \,.
$$
(b) The result (a) holds also for any projective $M$, that is without
the assumption of smoothness. \\
(c) Let $X$ and $\t X$ be two smooth projective manifolds and
$\pi:\t X\to X$ be a \ssc. Then
$$
HP\bu(X) \simeq HP\bu(\t X) \,.
$$
(d) Let $X,\t X,\pi$ be the same as in (c) and let $Z\subset X$ be an
arbitrary closed subset. Then
$$
\begin{array}{lcl}
HP\bu(Z)   & \simeq & HP\bu(\pi^{-1}(Z)) \,, \\
HP\bu(X,Z) & \simeq & HP\bu(\t X,\pi^{-1}(Z)) \,.
\end{array}
$$
\end{lemma}

\begin{proof}{}
We shall prove the propositions (a)--(d) by a simultaneous
induction in dimension. For $\dim E=0$ and $\dim X=1$ everything is obvious.
Suppose that (a)--(d) are proved when $\dim X<n$ and $\dim E<n-1$. Let us
prove these four propositions when $\dim E=n-1$ and $\dim X=\dim\t X=n$.

Consider a locally trivial fibration $\pi:E\to M$ where the fibers are all
isomorphic to the projective space $\C P^k$ for some $k\leqslant n-1$. Since
$\C P^k$ is birational to $(\C P^1)^{\times k}$  and by local triviality of
$\pi$ we conclude that $E$ is birational to the direct product $E':=M\times
(\C P^1)^{\times k}$. If $M$ is smooth as in part (a) of our statement, both
$E$ and $E'$ are smooth and the AKMW theorem (see Proposition \ref{A...W})
tells us that
$E$ and $E'$ can be related by a sequence of blow-ups and blow-downs. But
for $\dim E=\dim E'=n-1$, part (c) of the statement is applicable by our
induction hypothesis and we conclude that $HP\bu(E)=HP\bu(E')$. Finally,
$HP\bu(E')=HP\bu(M)$ according to Lemma \ref{cyl}. Thus, the
induction step is proved in part (a).

Let us now consider the fibration $\pi:E\to M$, $\dim E=n-1$, for an
arbitrary projective variety $M$ as in part (b). If $M$ is indeed singular
(perhaps even with intersecting components) we denote its
singular locus as $M_{\rm sing}\,$. By the Hironaka theorem there exists a
desingularization $\sigma:\t M\to M$, where $\t M$ consists of smooth
non-intersecting components and such that $M-M_{\rm sing}\simeq\t M-F$, where
$F:=\sigma^{-1}(M_{\rm sing})$. Let now $\t\pi:\t E\to \t M$ be the pull-back
of $\pi$ along $\sigma$. In this smooth situation, we have by proposition
(a) that $HP\bu(\t E)=HP\bu(\t M)$. Let us also consider the fibration
$\t\pi^{-1}(F)\to F$ (the restriction of $\t\pi$). Although its base $F$ may
be singular, its dimension ($\dim\t\pi^{-1}(F)<\dim E=n-1$) allows us to
use the induction hypothesis in part (b) to conclude that
$HP\bu(\t\pi^{-1}(F))=HP\bu(F)$. We want now to compare the polar homology
of the pair $\t M\supset F$ to that of $\t E\supset\t\pi^{-1}(F)$. The
isomorphisms $\pi_*: HP\bu(\t E)\simeq HP\bu(\t M)$ and $\pi_*:
HP\bu(\t\pi^{-1}(F))\simeq HP\bu(F)$ imply (as in Remark \ref{2pairs}) that
$$
HP\bu(\t E,\t\pi^{-1}(F)) \simeq HP\bu(\t M,F) \,.
$$
The varieties appearing in both sides of this equality have their birational
counterparts:
$$
\begin{array}{ccc}
\t M-F & \simeq & M-M_{\rm sing} \,, \\
\t E-\t\pi^{-1}(F) & \simeq & E-\pi^{-1}(M_{\rm sing}) \,.
\end{array}
$$
Hence, we can use Lemma \ref{BRE} to conclude that
\begin{equation} \label{*}
HP\bu(E,\pi^{-1}(M_{\rm sing})) \simeq HP\bu(M,M_{\rm sing})
\end{equation}
Since $\dim\pi^{-1}(M_{\rm sing})<\dim E=n-1$, we can apply the induction
hypothesis in part (b) to the fibration
$\pi^{-1}(M_{\rm sing})\to M_{\rm sing}$ and get the isomorphism
\begin{equation} \label{**}
HP\bu(\pi^{-1}(M_{\rm sing})) \simeq HP\bu(M_{\rm sing})
\end{equation}
Finally, the isomorphisms (\ref{*}) and (\ref{**}) and the map of
pairs $\pi:(E,\pi^{-1}(M_{\rm sing}))\to(M,M_{\rm sing})$ give the third
isomorphism $HP\bu(E)=HP\bu(M)$ as in Remark \ref{2pairs}, proving the
induction step in part (b).

Now, we turn to part (c) with two smooth projective varieties $X$ and
$\t X$, where $\dim X=\dim\t X=n$. It is sufficient to consider the case
when $\pi:\t X\to X$ is a single blow up with smooth center $M\subset X$.
Let us denote by $E=\pi^{-1}(M)\subset\t X$ the exceptional divisor.
Applying the proposition (a) to the fibration $\pi:E\to M$, we find that
$HP\bu(E)=HP\bu(M)$, while, by Lemma \ref{BRE}, we find that
$HP\bu(\t X,E)=HP\bu(X,M)$. These two isomorphisms imply the third one,
$HP\bu(\t X)=HP\bu(X)$, and we obtain the proof for part (c).

In part (d), we again consider the case of a single blowing up. Let $\pi$,
$X\supset M$, $\t X\supset E$ be the same as above and let $Z\subset X$ be
any closed subset. The subvariety $\pi^{-1}(Z)$ in $\t X$ may have many
components (even their dimensions may differ), so let us split these into
two groups, $\pi^{-1}(Z)=Z'\cup F$, where
$$
F=\pi^{-1}(Z\cap M) \,.
$$
In other words, $Z'$ is the union of the proper preimages of those components
of $Z$ not contained in $M$. So we have an isomorphism
$Z-Z\cap M\simeq Z'-Z'\cap F$, which by Lemma \ref{BRE} gives
$HP\bu(Z,Z\cap M)=HP\bu(Z',Z'\cap F)$. Besides, for
$\pi^{-1}(Z)=Z'\cup F$, we can write tautologically that
$HP\bu(Z',Z'\cap F)=HP\bu(\pi^{-1}(Z),F)$ and, hence,
$$
HP\bu(Z,Z\cap M)=HP\bu(\pi^{-1}(Z),F) \,.
$$
Taking into account that $HP\bu(F)=HP\bu(Z\cap M)$, which follows from (b)
for the fibration $F\to Z\cap M$, we conclude that
$$
HP\bu(Z)=HP\bu(\pi^{-1}(Z)) \,.
$$
The remaining equality, $HP\bu(X,Z)=HP\bu(\t X,\pi^{-1}(Z))$ follows
from (c), i.e.\ $HP\bu(X)=HP\bu(\t X)$, and by consideration of the map
of pairs $(\t X,\pi^{-1}(Z))\to (X,Z)$. Thus we have proved (d) and the whole
lemma.
\end{proof}

\begin{ssect}{}
If $V$ is a closed hypersurface in $X$, the embedding $i\!:V\hookrightarrow
X$ induces the corresponding homomorphisms in (co)homology. Namely, the
polar homology maps forward,
\begin{equation}\label{PHdirect}
i_*:HP_{q}(V) \to HP_{q}(X) \,.
\end{equation}
We have also the restriction map in sheaf cohomology,
$i^*:H^q(X,\OO_X)\to H^q(V,\OO_V)$.
If $V$ is smooth (or normal crossing), then
by Serre duality, $i^*$ produces the following
covariant homomorphism:
\begin{equation}\label{Hdirect}
i': H^{n-1-q}(V,K_V) \to H^{n-q}(X,K_X) \,.
\end{equation}
The proof of Theorem \ref{deRh} will be achieved essentially by comparing
the homomorphisms (\ref{PHdirect}) and (\ref{Hdirect}) and using (the simplest
case of) Lefschetz's hyperplane theorem. To describe this we begin
with a vanishing theorem.

\end{ssect}

\begin{prop} \label{Ko}
Let $V$ be an ample divisor and $D$ be a normal crossing divisor in a smooth
projective manifold $X$. Then
$$
H^p(X,K_X(V+D))=0\,, ~p>0\,.
$$
\end{prop}

\noindent This mild generalization (i.e. to $D\ne\varnothing$) of the Kodaira
vanishing theorem can be found in ref.
\cite{some}. Now suppose also that $V$ is a normal crossing
divisor. Then the long exact sequence in cohomology of
\begin{equation} \label{KVD}
0 \to K_X(D) \to K_X(V+D) \to K_V(D) \to 0  \,,
\end{equation}
gives the following.

\begin{prop} \label{Lefsch1}
If $V$ and $D$ are normal crossing divisors
in a smooth projective $X$, with $V$ ample, then
$$
\begin{array}{rlclr}
i': & H^p(V,K_V(D)) & \toiso & H^{p+1}(X,K_X(D)) & {\rm for}~ p>0 \,, \\
i': & H^0(V,K_V(D)) & \twoheadrightarrow & H^1(X,K_X(D)) \,. &
\end{array}
$$
\end{prop}

\begin{prop} \label{Lefsch2}
If $V$ is an ample normal crossing subvariety
in a smooth projective $X$ and $m={\rm codim}\, V$, then
$$
\begin{array}{rlclr}
i': & H^p(V,K_V) & \toiso & H^{p+m}(X,K_X) & {\rm for}~ p>0 \,, \\
i': & H^0(V,K_V) & \twoheadrightarrow & H^m(X,K_X) \,. &
\end{array}
$$
\end{prop}

\noindent
This follows trivially from the Lefschetz theorem (Proposition
\ref{Lefsch1}) by considering a flag
$V=V^m\subset V^{m-1}\subset\ldots\subset V^1\subset V^0=X$
with $V^{i+1}$ being an ample normal crossing divisor in $V^i$ (such a flag
exists by definition).

\begin{prop} \label{Lefsch3}
Let $V=V^m\subset V^{m-1}\subset\ldots\subset V^1\subset V^0=X$
be as above and let $D\subset X$ be a normal crossing divisor which
intersects each $V^i$ transversely (so that $D\cap V^i$ is also a
normal crossing divisor in $V^i$). Then
$$
\begin{array}{rlclr}
i': & H^p(V,K_V(D)) & \toiso & H^{p+m}(X,K_X(D)) & {\rm for}~ p>0 \,, \\
i': & H^0(V,K_V(D)) & \twoheadrightarrow & H^m(X,K_X(D)) \,. &
\end{array}
$$
\end{prop}

\begin{rem}
Suppose Theorem \ref{rel} is proven. Then Proposition \ref{Lefsch1}
has also a similar implication in polar homology
(with $D=\varnothing$), namely:
$$
\begin{array}{rlclr}
i_*: & HP_k(V) & \toiso & HP_k(X) & {\rm for}~ k < n-1 \,, \\
i_*: & HP_{n-1}(V) & \twoheadrightarrow & HP_{n-1}(X) \,. &
\end{array}
$$ It may be interesting to notice that this has the following
topological analogue. For an $n$-dimensional $CW$-complex $X$ and its
$(n-1)$-skeleton $i\!:V \hookrightarrow X$, the map $i_*:H_{q}(V)\to
H_{q}(X)$ is an isomorphism of cellular homology for $0\leqslant q < n-1$
and is surjective for $q=n-1$.

Thus, by Lefschetz's theorem in the form of Proposition \ref{Lefsch1} one
can view an {\it ample divisor} in the context of polar homology as an
analogue of the $(n-1)$-{\it skeleton} in topology. Of course, the Morse
theory proof of the Lefschetz theorem shows that the topological $(n-1)$-skeleton
can indeed be taken to lie in the hyperplane.

\end{rem}

\begin{ssect}{Proof of Theorem \ref{deRh}.} \label{proof}
Let us show first that the map $\rho$ in eq.\ (\ref{rho}) is surjective.
Take an arbitrary ample smooth subvariety
$i\!:V\hookrightarrow X$,
$\dim V=q$. Then
$i'\!:H^0(V,K_V)\twoheadrightarrow H^{n-q}(X,K_X)$ is surjective by the
Lefschetz theorem \ref{Lefsch2}.
But each element $\alpha\in H^0(V,K_V)$ corresponds, by definition,
to a cycle $a=(V,i,\alpha)$ in $HP_q(X)$ and $\rho([a])=i'(\alpha)$.
Thus $\rho$ is onto.

To prove injectivity we must show that for a $q$-cycle $a$ the vanishing
$\rho([a])=0\in H^{n-q}(X,K_X)$ implies that $a=\partial b$ for some polar
$(q+1)$-chain $b$. Let $a=\sum_k(A_k,f_k,\alpha_k)\in{\cal C}_q(X),\:
\partial a=0$, be an arbitrary $q$-cycle. Its support,
${\rm supp}\,a=Z=\cup_kZ_k$, may be a singular reducible subvariety\footnote{
We may suppose without loss of generality that $Z$ has the same number of
components as the number of terms in $a=\sum_k(A_k,f_k,\alpha_k)$.}
in $X$. Let $Z_{\rm sing}$ be the subset of singular points of $Z$
(including, of course, possible points of intersection of its components).
By the Hironaka theorem we can find a
blow-up $\pi:\t X\to X$ such that the following conditions are satisfied.
\begin{description}

\item a) There is a $q$-dimensional subvariety $\t Z\subset\t X$ which consists of
smooth non-intersecting  components and such that $\pi(\t Z)=Z$ and $\pi$
gives us a birational map of $\t Z$ onto $Z$.

\item b) $\t Z$ is included into a nested sequence of subvarieties:
\begin{equation} \label{ZY-flag}
\t Z\subset\t Y=V^{n-q-1}\subset V^{n-q-2}\subset\ldots\subset V^1\subset V^0=\t X
\end{equation}
where ${\rm codim}\,V^i=i$ (in particular, $\dim\t Y=q+1$) and each
$V^{i+1}$ is an ample normal crossing divisor in $V^i$, so that $\t Y$, in
particular, is an ample \ncsv\ in $\t X$. (If $q=n$ our proposition is
obvious: $HP_n(X)=H^0(X,K_X)$, while for $q=n-1$ we set simply $\t Y=\t X$.)

\item c) The preimage $D:=\pi^{-1}(Z_{\rm sing})$ of the singular locus of $Z$ is
a normal crossing divisor in $\t X$ which also intersects transversely $\t
Z$, $\t Y$ as well as all other elements $V^i$ of the flag (\ref{ZY-flag}).

\end{description}
We can ensure this by applying the Proposition \ref{Hi-flag} to each
component of $Z$. The possibility to satisfy the condition (c) is also
guaranteed by the Hironaka theorem. After that we can achieve the
ampleness of $V^1,V^2,\ldots,\t Y$ by adding
sufficiently ample components to them, which can be done
preserving normal crossings.

We are now prepared to replace the original polar cycle $a\in{\cal C}_q(X)$,
which has a singular support $Z\subset X$, with a cycle supported on $\t Z$
in $\t X$. Recall that $\t Z$ may have several components, $\t Z=\cup_k\t
Z_k$, but these do not intersect. Each $q$-dimensional smooth subvariety
$i_k:\t Z_k\hookrightarrow\t X$ acquires a meromorphic $q$-form $\t\alpha_k$
defined on $\t Z_k$. This can be seen
by noticing that there exists a smooth manifold $\t A_k$ birational to $A_k$
with a commutative square
$$
\begin{CD}
\t A_k    @>>>    \t Z_k      \\
  @VVV              @V{\pi}VV   \\
   A_k    @>{f_k}>>       Z_k
\end{CD}
$$
which allows us to pull back $\alpha_k$ from $A_k$ to $\t A_k$ and then to
push it forward to $\t Z_k$. We claim that each triple $(\t Z_k,\t i_k,\t\alpha_k)$
is admissible. Since $a$ was a closed chain, the polar locus of
$\alpha_k$ was mapped by $f_k$ to $Z_{\rm sing}$. Therefore, we have that
$\dvsr_\infty\t\alpha_k\subset \t Z_k\cap D$, where
$D=\pi^{-1}(Z_{\rm sing})$. By virtue of c) above, this guarantees that
the polar divisor is normal crossings. Thus we need now only show that
$\t\alpha_k$ has at most first order poles. The form $\t\alpha_k$ is obtained
from $\alpha_k$ by means of pushforwards and pullbacks, which we claim
both preserve the property of having only first order poles. The first follows
from a local calculation with the cover $z\mapsto z^n$ about
the smooth locus of a branch divisor. For the second we use
the observation that for top degree forms, having first
order poles is the same as being logarithmic, where logarithmic forms are
locally generated as a ring by forms $df/f=d\log f$ and so are also logarithmic
on pullback.

So each $(\t Z_k,\t i_k,\t\alpha_k)$ defines a $q$-chain in $\t X$. However,
the sum of these
triples, $\t a=\sum_k(\t Z_k,\t i_k,\t\alpha_k)$, does not necessarily form a
cycle\footnote{
For example, a 1-cycle in $X$ can be supported on a self-intersecting
rational  curve $Z$. Then the resolved smooth curve $\t Z\subset\t X$ will
be equipped with a  meromorphic 1-form which has simple poles at
the resolution of the double point of $Z$ and, hence, the resolved curve is
no longer a cycle in $\t X$.}.
Nevertheless, $\t a$ has no boundary modulo $D$ in $\t X$, so
we consider $\t a$ as a $q$-cycle in ${\cal C}_q(\t X,D)$.

Now we suppose that $\rho([a])=0\in H^{n-q}(X,K_X)$ and try to prove
that $[a]=0$ in $HP_q(X)$. Let us note first that by (\ref{l.seq.rel.})
it is enough to prove the vanishing of $[a]$ modulo
$Z_{\rm sing}\subset X$, that is in $HP_q(X,Z_{\rm sing})$, because
$\dim Z_{\rm sing}<\dim Z=q$ and so $HP_q(Z_{\rm sing})=0$. Secondly, since
$\pi_*:HP_q(\t X,D)\toiso HP_q(X,Z_{\rm sing})$
by Lemma \ref{bl-up}(d) and since, obviously, $\pi_*[\t a]=[a]$ it is
sufficient to prove that $[\t a]=0\in HP_q(\t X,D)$. To prove this latter
vanishing we have to show first that $\t\rho([\t a])=0$, where
$$
\t\rho: HP_q(\t X,D)\to H^{n-q}(\t X,K_{\t X}(D))
$$
is the relative analogue of the map $\rho$, which is the subject of the
proposition under consideration (cf.\ eqs.\ (\ref{rho}) and (\ref{rho
rel})). For this aim let us collect the relevant maps in polar homology
recalling the isomorphisms in Lemma \ref{bl-up} as well as the isomorphism
$H^{n-q}(\t X, K_{\t X})\toiso H^{n-q}(X, K_X)$, which holds for smooth
birationally equivalent $X$ and $\t X$, in the following commutative diagram:
\cd
Then, from $\rho([a])=0$, it follows that
$\t\rho([\t a])=0\in H^{n-q}(\t X, K_{\t X}(D))$.

We are ready now to finish the proof. To simplify the notations let us write
$\t a=(\t Z,\t i,\t\alpha)$ for the sum $\sum_k(\t Z_k,\t i_k,\t\alpha_k)$,
where $\t\alpha\in H^0(\t Z,K_{\t Z}(D))$, while
$\t i\!:\t Z\hookrightarrow\t X$ is the embedding of the union of smooth
non-intersecting components $\t Z=\cup_k\t Z_k$ into $\t X$. The map
$\t\rho$ applied to $\t a$ corresponds to the map
$\t i'\!: H^0(\t Z,K_{\t Z}(D)) \to H^{n-q}(\t X, K_{\t X}(D))$,
that is to say, $\t\rho([\t a])=\t i'(\t\alpha)$. Thus, we have that
$\t i'(\t\alpha)=0$.
Since the embedding $\t i\!:\t Z\hookrightarrow\t X$ can be described as a
composition of two embeddings, $\t i_{\t Y}:\t Z\hookrightarrow\t Y$ and
$\t j\!:\t Y\hookrightarrow\t X$ the above map $\t i'$
factors in this case through $H^1(\t Y,K_{\t Y}(D))$:
\begin{equation}\label{seq1}
\begin{CD}
H^0(\t Z,K_{\t Z}(D))          @>\t i'_{\t Y}>>
H^1(\t Y,K_{\t Y}(D))       @>\t j'>\mbox{\large$\sim$}>
H^{n-k}(\t X,K_{\t X}(D))  \,,
\end{CD}
\end{equation}
where $\t j'\circ\t i'_{\t Y}=\t i'$ and $\t j'$ is an isomorphism by the
ampleness of $\t Y$ (see Proposition \ref{Lefsch3}). It follows that
$\t i'_{\t Y}(\t\alpha)=0$ and the problem reduces to a codimension one
situation: $\t Z\subset\t Y$. We can consider now the following exact sequence:
\begin{equation}\label{seq2}
\begin{CD}
0\to K_{\t Y}(D)\to K_{\t Y}(D\cap\t Y+\t Z)
@>res_{\t Z}>>
K_{\t Z}(D) \to 0
\end{CD}
\end{equation}
and the corresponding long sequence in cohomology. The latter allows us
to conclude that the vanishing
$\t i'_{\t Y}(\t\alpha)=0,~ \t\alpha\in H^0(\t Z,K_{\t Z}(D))$,
implies that $\t\alpha=\res_{\t Z}\,\t\beta$ for some
$\t\beta\in H^0(\t Y,K_{\t Y}(D\cap\t Y+\t Z))$. In terms of polar chains in
$\t X$ (modulo $D$), this means that
$\t a=(\t Z,\t i,\t\alpha)=\partial(\t Y,\t j,\t\beta)$, or
$[\t a]=0\in HP_q(\t X,D)$. As we explained above, this implies that
$[a]=0\in HP_q(X)$, which proves the injectivity of $\rho$.
$\blacksquare$

\end{ssect}



\bigskip

\begin{ackn}
We are mostly indebted to D.~Orlov for very detailed discussions and
especially for {\it disproving} an assertion, which was a technical, but
extremely misleading point on our way. We would like to thank also
A.~Kuznetsov and A.~Pukhlikov for discussions of related problems, and the
referee for a number of suggestions.

A.R.\ and B.K.\ are grateful  to the MPI f\"ur Mathematik in Bonn,
as well as to the ESI in Vienna and  IHES in Bures-sur-Yvette, for
kind hospitality during the work on this paper.

The work of B.K.\ was partially supported by PREA of Ontario, McLean and
NSERC research grants.
The work of A.R.\ was supported in part by the Grants RFBR-01-01-00539,
INTAS-99-1705 and the Grant 00-15-96557 for the support of scientific
schools.
The work of R.T.\ was supported by a Royal Society university research
fellowship.

\end{ackn}


\begin{smallbibl}{AKMW}

\bibitem[AKMW]{A...W} D.~Abramovich, K.~Karu, K.~Matsuki,
and J.~W\l odarczyk,
{\em Torification and factorization of birational maps}\/,
J.\ Amer.\ Math.\ Soc.\ {\bf 15} (2002) 531--572
[math.AG/9904135]

\bibitem[DT]{DT} S.K.~Donaldson and R.P.~Thomas,
{\em Gauge theory in higher dimensions}\/,
In: The Geometric Universe:
Science, Geometry and the work of Roger Penrose, S. A. Huggett et al
(eds), Oxford University Press, Oxford (1998)

\bibitem[EV]{some} H.~Esnault and E.~Viehweg,
{\em Lectures on vanishing theorems}\/,
DMV Seminar, {\bf 20}, Birkhauser Verlag, Basel, (1992)

\bibitem[G]{Griff} P.A.~Griffiths, {\em Variations on a Theorem of
Abel}\/, Invent.\ Math., {\bf 35} (1976) 321--390

\bibitem[H]{Hi} H.~Hironaka,
{\em Resolution of singularities of an algebraic variety
over a field of characteristic zero, I}\/,
Ann.\ of Math.\ (2) {\bf 79} (1964) 109--203; \\
H.~Hironaka,
{\em Resolution of singularities of an algebraic variety
over a field of characteristic zero, II}\/,
Ann.\ of Math.\ (2) {\bf 79} (1964) 205--326

\bibitem[KR1]{KR}  B.~Khesin and A.~Rosly, {\em Polar homology}\/,
[math.AG/0009015]

\bibitem[KR2]{KR2} B.~Khesin and A.~Rosly,
{\em Polar Homology and Holomorphic Bundles,} \\
Phil.\ Trans.\ R.\ Soc.\ Lond., Ser.\ A, {\bf 359}, (2001) 1413--1428
[math.AG/0102152]

\bibitem[W]{W}  J.~W\l odarczyk,
{\em Combinatorial structures on toroidal varieties and a proof of the weak
Factorization Theorem}\/,
[math.AG/9904076]

\end{smallbibl}

\end{document}